# GENERALIZED URN MODELS OF EVOLUTIONARY PROCESSES[1]

By Michel Benaïm, Sebastian J. Schreiber and Pierre Tarrès

*Institut de Mathématiques, College of William and Mary, and Université Paul Sabatier*

Generalized Pólya urn models can describe the dynamics of finite populations of interacting genotypes. Three basic questions these models can address are: Under what conditions does a population exhibit growth? On the event of growth, at what rate does the population increase? What is the long-term behavior of the distribution of genotypes? To address these questions, we associate a mean limit ordinary differential equation (ODE) with the urn model. Previously, it has been shown that on the event of population growth, the limiting distribution of genotypes is a connected internally chain recurrent set for the mean limit ODE. To determine when growth and convergence occurs with positive probability, we prove two results. First, if the mean limit ODE has an "attainable" attractor at which growth is expected, then growth and convergence toward this attractor occurs with positive probability. Second, the population distribution almost surely does not converge to sets where growth is not expected and almost surely does not converge to "nondegenerate" unstable equilibria or periodic orbits of the mean limit ODE. Applications to stochastic analogs of the replicator equations and fertility-selection equations of population genetics are given.

**1. Introduction.** The founder-effect in population genetics refers to the establishment of a new population consisting of a few founders that is isolated from the original population. A founder-effect can occur when a small number of individuals colonize a place previously uninhabited by their species. In this case, the founding population is geographically isolated from the original population. A founder-effect due to temporal isolation can occur when a population passes through a bottle neck after which only a few

Received September 2001; revised May 2003.
[1]Supported in part by NSF Grant DMS-00-77986.

individuals survive. Several fundamental questions surrounding the founder-effect include: What is the probability that a founding population successfully establishes itself? If a founding population establishes itself, what is the population's growth rate and what is the long-term genotypic or phenotypic composition of the population? How does the initial genotypic composition and initial population size influence the likelihoods of the various outcomes?

To address these questions, we consider a general class of Pólya urn models introduced in [12]. Traditionally, Pólya urn models are described as involving an urn, which contains a finite number of balls of different colors. At discrete moments in time, balls are added or removed from the urn according to probabilities that only depend on their distribution and number at that point in time. The pertinence of these models to evolutionary questions is self-evident if we view the balls as individuals whose color represents their genotype or phenotype, adding or removing balls as replication or death of individuals and updates of the urn as interactions between individuals. When balls are added at a constant rate (i.e., a fixed number of individuals are added at every update of the process), Pólya urn models have been studied extensively by Arthur, Ermol'ev and Kaniovskii [1], Benaïm and Hirsch [5, 6], Hill, Lane and Sudderth [8], Pemantle [10] and Posch [11] amongst others. What makes the models introduced in [12] more relevant to population processes is that they permit the removal of balls, as well as the addition of balls at nonconstant rates. Consequently, extinction of the entire population or one or more subpopulations is possible in finite time.

The article is structured as follows. In Section 2 we introduce the class of generalized urn models. As examples, replicator processes and fertility-selection processes with and without mutations are introduced. These urn models typically predict that either the populations go extinct or grow, and that demographic stochasticity is most pronounced when the population is small. Once the population starts to get large, it tends to grow in an essentially deterministic fashion. For this reason, on the event of nonextinction, the dynamics of the distribution of types in the population (i.e., the distribution of the color of balls in the urn) are strongly correlated to the asymptotic behavior of an appropriately chosen ordinary differential equation (ODE), commonly called the *mean limit ODE*. In Section 3 we describe the mean limit ODE, and recall a theorem [12] that on the event of growth relates the asymptotic behaviors of the stochastic process to its mean limit ODE. Using this result, we derive a time averaging principle and a competitive exclusion principle for replicator processes. We also show that additive fertility-selection processes almost surely converge on the event of nonextinction to a fixed point of the mean limit ODE. These results, however, provide no insight into when population growth occurs with positive probability and which limiting behaviors occur with positive probability. In [12], Theorem 2.6, it was shown that if the population *grows with probability* 1 (i.e., more balls are being added than removed at each update), then the population distribution converges with positive probability to "attainable" attractors for the mean limit ODE. However, for most evolutionary and ecological processes extinction occurs with positive probability and, consequently, the almost surely growth assumption is not meet. In Section 4 we remedy this issue and prove that if there is an "attainable" attractor for the mean limit ODE *at which growth is expected* (i.e., on average more balls are added than removed), then there is a positive (typically less than 1) probability that the population grows and the population distribution converges to the attractor. In addition, we provide an estimate for this probability and for the rate of growth. In Section 5 we prove nonconvergence to invariant subsets of the mean limit ODE where growth is not expected and nonconvergence to "nondegenerate" unstable equilibria and periodic orbits of the mean limit ODE. In Section 6 we combine our results to give necessary and sufficient conditions for growth with positive probability for processes with gradient-like mean limit ODEs and apply these conditions to additive fertility processes with mutation.

**2. Generalized urn models.** In this section we introduce a class of generalized urn models. This class generalizes the urn models introduced in [12]. Three examples or evolutionary processes described by these generalized urn models are given in the sections below. One of these processes, replicator processes, was described previously in [12] and stochastic versions of the replicator equations [9]. The other two processes, fertility-selection processes with and without mutation, are stochastic versions of the fertility-selection equations with and without mutation [9].

Due to the fact that we are dealing with finite populations consisting of individuals that are one of $k$ types, we consider Markov chains on the positive cone

$$\mathbb{Z}_+^k = \{z = (z^1, \ldots, z^k) \in \mathbb{Z}^k : z^i \geq 0 \text{ for all } i\}$$

of the set $\mathbb{Z}^k$ of $k$-tuples of integers. Given a vector $w = (w^1, \ldots, w^k) \in \mathbb{Z}^k$, define

$$|w| = |w^1| + \cdots + |w^k| \quad \text{and} \quad \alpha(w) = w^1 + \cdots + w^k.$$

We shall always write $\|\cdot\|$ for the Euclidean norm on $\mathbb{R}^k$.

Let $z_n = (z_n^1, \ldots, z_n^k)$ be a homogeneous Markov chain with state space $\mathbb{Z}_+^k$. In our context, $z_n^i$ corresponds to the number of balls of color $i$ at the $n$th update. Associated with $z_n$ is the random process $x_n$ defined by

$$x_n = \begin{cases} \dfrac{z_n}{|z_n|}, & \text{if } z_n \neq 0, \\ 0, & \text{if } z_n = 0, \end{cases}$$

which is the distribution of balls at the $n$th update. Note that when there are no balls at the $n$th update, we set $x_n$ to zero which we view as the "null" distribution. Let $S_k \subset \mathbb{R}^k$ denote the unit $k-1$ simplex, that is,

$$S_k = \left\{ x = (x^1, \ldots, x^k) \in \mathbb{R}^k : x^i \geq 0, \sum_{i=1}^k x^i = 1 \right\}.$$

Let $\Pi : \mathbb{Z}_+^k \times \mathbb{Z}_+^k \mapsto [0,1]$ denote the transition kernel of the Markov chain $z_n$. In other words, $\Pi(z, z') = P[z_{n+1} = z' | z_n = z]$. We place the following assumptions on the Markov chains $z_n$ of interest:

(A1) At each update, there is a maximal number of balls that can be added or removed. In other words, there exists a positive integer $m$ such that $|z_{n+1} - z_n| \leq m$ for all $n$.

(A2) There exist Lipschitz maps

$$\{p_w : S_k \to [0,1] : w \in \mathbb{Z}^k, |w| \leq m\}$$

and a real number $a > 0$, such that

$$|p_w(z/|z|) - \Pi(z, z+w)| \leq a/|z|$$

for all nonzero $z \in \mathbb{Z}_+^k$ and $w \in \mathbb{Z}^k$ with $|w| \leq m$.

Since we view updates of the Markov chain to correspond to the effect of interactions between individuals, assumption (A1) implies that each interaction results in the addition or removal of no more than a maximum number of individuals. Assumption (A2) assures that there is a well-defined mean limit ODE for the urn models.

2.1. *Replicator processes.* Consider a system consisting of a finite population of individuals playing $k$ different strategies. At each update of the population, pairs of individuals are chosen randomly with replacement from the population. The chosen individuals replicate and die according to probabilities that only depend on their strategies. More precisely, let $m$ be a nonnegative integer that represents the maximum number of progeny that any individual can produce in one update. Let $\{R_n\}_{n \geq 0}$ and $\{\tilde{R}_n\}_{n \geq 0}$ be sequences of independent identically distributed random $k \times k$ matrices whose entries take values in the set $\{-1, 0, 1, \ldots, m\}$. Let $R_n^{ij}$ and $\tilde{R}_n^{ij}$ denote the $ij$th entry of $R_n$ and $\tilde{R}_n$, respectively. Let $\{r_n\}_{n \geq 0}$ be a sequence of independent identically distributed random $k \times 1$ matrices whose entries take values in the set $\{-1, 0, 1, \ldots, m\}$. We define a replicator process according to the following rules:

1. Two individuals are chosen at random with replacement from the the population. Make note of the individuals chosen and return them to the population.

2. If the same individual is chosen twice and it plays strategy $i$, then $r_n^i$ individuals of strategy $i$ are added to population.
3. If two distinct individuals are chosen, say strategy $i$ and strategy $j$, then add $R_n^{ij}$ individuals of strategy $i$ and add $\tilde{R}_n^{ji}$ individuals of strategy $j$.

It is not difficult to verify that this process satisfies assumptions (A1) and (A2) with $p_w(x) = x^i x^i P[R_1^{ii} + \tilde{R}_1^{ii} = w^i] + \sum_{j \neq i} 2x^i x^j P[R_1^{ij} = w^i] P[\tilde{R}_1^{ji} = 0]$, whenever $w = (0, \ldots, 0, w^i, 0, \ldots, 0)$ with $w^i \neq 0$, $p_w(x) = 2x^i x^j P[R_1^{ij} = w^i] P[\tilde{R}_1^{ji} = w^j]$, whenever $w = (0, \ldots, 0, w^i, 0, \ldots, 0, w^j, 0, \ldots, 0)$ with $w^i \neq 0$, $w^j \neq 0$, and $i \neq j$, $p_w(x) = 0$, whenever $w$ has three or more nonzero coordinates, and $p_0(x) = 1 - \sum_{w \neq 0} p_w(x)$.

2.2. *Fertility-selection processes.* Consider a population of diploid individuals that has $k$ distinct alleles $A_1, \ldots, A_k$ that occupy a single locus. Assume that the population is monoecious (i.e., there is only one "sex") and that each individual chooses its mate randomly from the population. We assume that individuals die immediately after mating. Although there is no distinction between individuals of genotype $A_i A_j$ and $A_j A_i$, we develop an urn model of the form $z_n = (z_n^{ij}) \in \mathbb{Z}_+^{k \times k}$ as it is notationally more convenient. For $i \neq j$, let $z_n^{ij} = z_n^{ji}$ denote the number of individuals of genotype $A_i A_j$ at the $n$th update of the population. Alternatively, let $z_n^{ii}$ denote twice the number of individuals of genotype $A_i A_i$ at the $n$th update of the population. Hence, $|z_n| = \sum_{i,j} z_n^{ij}$ equals twice the total number of individuals in the population.

For every pair of genotypes, say $A_i A_j$ and $A_r A_s$, we associate a sequence of i.i.d. random variables $G_n(ij, rs) = G_n(rs, ij)$ that take values in $\{0, \ldots, m\}$, where $m$ represents the maximal number progeny produced by a mating and where $G_n(ij, rs)$ represents the number of progeny produced by a mating between genotypes $A_i A_j$ and $A_r A_s$ at update $n$. Let $z_n \in \mathbb{Z}_+^{k \times k}$ be a Markov chain satisfying $z_n^{ij} = z_n^{ji}$ and updated according to the following rules:

1. If the population size is less than two, then the population goes extinct. In other words, if $|z_n| < 4$, then $z_{n+1} = 0$.
2. Pick two individuals at random without replacement from the population, say genotypes $A_i A_j$ and $A_r A_s$.
3. Remove the chosen individuals from the population (i.e., they die).
4. Add $G_n(ij, rs)$ individuals to the population. The genotype, say $A_u A_v$, of each added individual is independently determined by random mating probabilities (i.e., $u$ equals $i$ or $j$ with equal probability and $v$ equals $r$ or $s$ with equal probability).

Define
$$x_n^{ij} = \begin{cases} \dfrac{z_n^{ij}}{|z_n|}, & \text{if } z_n \neq 0, \\ 0, & \text{if } z_n = 0. \end{cases}$$

Hence, if $i \neq j$, then $2x^{ij} = 2x^{ji}$ equals the proportion of the population with genotype $A_i A_j$. Alternatively, $x^{ii}$ is the proportion of the population with genotype $A_i A_i$.

To see that this process satisfies assumption (A2), let $w$ and $z$ be in $\mathbb{Z}_+^{k \times k}$ such that $w^{ij} = w^{ji}$ and $z^{ij} = z^{ji}$ for all $1 \leq i, j \leq k$. Assume that $|z| > 0$. $p_w(z/|z|)$ is given by a linear combination of the terms $x^{ij} x^{rs}$. On the other hand, $\Pi(z, z+w)$ is given by the corresponding linear combination of the terms $x^{ij} x^{rs} |z|/(|z|-2)$, when $\{i,j\} \neq \{r,s\}$ and $x^{ij}(x^{rs}|z|-2)/(|z|-2)$. From these observations it follows that assumption (A2) is satisfied.

2.3. *Fertility-selection process with mutations.* To account for mutations in the fertility-selection process, let $\mu(ij, rs) \geq 0$ for $1 \leq i, j, r, s \leq k$ be such that $\sum_{r \leq s} \mu(ij, rs) = 1$ for all $1 \leq i \leq j \leq k$. The quantity $\mu(ij, rs)$ represents the probability that the genotype $A_i A_j$ mutates to the genotype $A_r A_s$. The fertility process with mutation is given by the first three rules of the fertility process without mutation and replacing the fourth rule with: Add $G_n(ij, rs)$ individuals to the population. For each added individual, its genotype is determined by two steps. First, determine a genotype $A_u A_v$ according to random mating probabilities. The probability that the added individual has genotype $A_{\tilde u} A_{\tilde v}$ is given by $\mu(uv, \tilde u \tilde v)$. For reasons similar to the fertility-selection process without mutation, this process also satisfies assumptions (A1) and (A2).

**3. Mean limit ODEs.** To understand the limiting behavior of the $x_n$, we express $x_n$ as a stochastic algorithm using the following lemma.

LEMMA 1. *Let $z_n$ be a Markov chain on $\mathbb{Z}_+^k$ satisfying assumptions* (A1) *and* (A2) *with mean limit transition probabilities $p_w : S_k \to [0, 1]$. Let $\mathcal{F}_n$ denote the $\sigma$ field generated by $\{z_0, z_1, \ldots, z_n\}$. There exists sequences of random variables $\{U_n\}$ and $\{b_n\}$ adapted to $\mathcal{F}_n$, and a real number $K > 0$ such that:*

(i) *if $z_n \neq 0$, then*

$$(1) \qquad x_{n+1} - x_n = \frac{1}{|z_n|} \left( \sum_{w \in \mathbb{Z}^k} p_w(x_n)(w - x_n \alpha(w)) + U_{n+1} + b_{n+1} \right),$$

(ii) $E[U_{n+1}|z_n] = 0$,
(iii) $\|U_n\| \leq 4m$ and $E[\|U_{n+1}\|^2 | \mathcal{F}_n] \leq 4m^2$,
(iv) $\|b_{n+1}\| \leq \frac{K}{\max\{1, |z_n|\}}$.

PROOF. The proof of this lemma is very similar to the proof of Lemma 2.1 in [12]. Consequently, we only provide an outline of the proof. Define

$$g(x) = \sum_{w \in \mathbb{Z}^k} p_w(x)(w - x\alpha(w)),$$

$$U_{n+1} = (x_{n+1} - x_n - E[x_{n+1} - x_n|z_n])|z_n|,$$
$$b_{n+1} = |z_n|E[x_{n+1} - x_n|z_n] - g(x_n).$$

From these definitions it follows that (i) and (ii) hold. For the remainder of the proof, let $z = z_n$ and $x = x_n$. To prove (iii), notice that if $z_{n+1} \neq 0$, then it can be shown that

(2) $$\|(x_{n+1} - x)|z|\| \leq 2m,$$

as no more than $m$ balls are being added or removed at any update. Alternatively, if $z_{n+1} = 0$, then it must be that $|z_n| \leq m$ since no more than $m$ balls can be removed at a single update. In which case, $x_{n+1} = 0$ and $\|(x_{n+1} - x_n)|z_n|\| = \|z_n\| \leq m$. From this (ii) follows. To prove (iv), notice that if $z \neq 0$, then

(3) 
$$|z|E[x_{n+1} - x|z_n = z]$$
$$= \sum_{w \neq -z} \frac{|z|(w - \alpha(w)x)}{|z| + \alpha(w)} \Pi(z, z+w) - z\Pi(z, 0).$$

If $|z| > m$, (3) can be used to show

$$\|b_{n+1}\| \leq \left\| \sum_{|w| \leq m} \left( \frac{|z|}{|z| + \alpha(w)} \Pi(z, z+w) - p_w(x) \right)(w - \alpha(w)x) \right\|.$$

Applying assumptions (A1) and (A2) implies that there is $K_1 > 0$ such that $\|b_{n+1}\| \leq K_1/|z_n|$, whenever $|z_n| > m$. On the other hand, if $|z| \leq m$, then the definition of $b_{n+1}$ implies that $\|b_{n+1}\| \leq 2m + \sup_{x \in S_n} \|g(x)\|$. Choosing $K$ sufficiently larger than $K_1$ completes the proof of (iv). □

The recurrence relationship (1) can be viewed as a "noisy" Cauchy–Euler approximation scheme with step size $1/|z_n|$ for solving the ordinary differential equation

(4) $$\frac{dx}{dt} = \sum_{w \in \mathbb{Z}^k} p_w(x)(w - x\alpha(w)),$$

which we call the *mean limit ODE*. When the number of individuals in the population grow without bound, the step size decreases to zero and it seems reasonable that there is a strong relationship between the limiting behavior of the mean limit ODE and the distribution of balls $x_n$. To make the relationship between the stochastic process $x_n$ and the mean limit ODE more transparent, it is useful to define a continuous time version of $x_n$ where time is scaled in an appropriate manner. Since the number of events (updates) that occur in a given time interval is likely to be proportional to the size of the population, we define *the time $\tau_n$ that has elapsed by update $n$* as

$$\tau_0 = 0,$$

$$\tau_{n+1} = \begin{cases} \tau_n + \dfrac{1}{|z_n|}, & \text{if } z_n \neq 0, \\ \tau_n + 1, & \text{if } z_n = 0. \end{cases}$$

The continuous time version of $x_n$ is given by

(5) $$X_t = x_n \quad \text{for } \tau_n \leq t < \tau_{n+1}.$$

To relate the limiting behavior of the flow $\phi_t(x)$ of (4) to the limiting behavior of $X_t$, Schreiber [12] proved the next theorem using techniques of Benaïm [2, 4]. Recall, a set $C$ is called *invariant* for the flow $\phi_t$ provided that $\phi_t(C) = C$ for all $t \in \mathbb{R}$. A compact invariant set $C$ is *internally chain recurrent* provided that for every $x \in C$, $T > 0$ and $\varepsilon > 0$, there exist points $x_1, x_2 \ldots, x_s$ in $C$ and times $t_1, \ldots, t_s$ greater than $T$ such that $x_1 = x_s = x$ and $\|\phi_{t_i}(x_i) - x_{i+1}\| < \varepsilon$ for $1 \leq i \leq s - 1$. Given a function $X_t : \mathbb{R}_+ \to \mathbb{R}^k$ or a sequence $\{x_n\}_{n \geq 0}$ in $\mathbb{R}^k$, we define the *limit sets*, $L(X_t)$ and $L(x_n)$, of $X_t$ and $x_n$. $L(X_t)$ is the set of $p \in \mathbb{R}^k$ such that $\lim_{k \to \infty} X_{t_k} = p$ for some subsequence $\{t_k\}_{k \geq 0}$ with $\lim_{k \to \infty} t_k = \infty$. $L(x_n)$ is the set of $p \in \mathbb{R}^k$ such that $\lim_{k \to \infty} x_{n_k} = p$ for some subsequence $\{n_k\}_{k \geq 0}$ with $\lim_{k \to \infty} n_k = \infty$.

THEOREM 1 ([12]). *Let $z_n$ be a Markov process satisfying assumptions* (A1) *and* (A2) *with mean limit ODE* (4). *Then on the event* $\{\liminf_{n \to \infty} \frac{|z_n|}{n} > 0\}$:

1. *The interpolated process $X_t$ is almost surely an* asymptotic pseudotrajectory *for the flow $\phi_t$ of the mean limit ODE. In other words, $X_t$ almost surely satisfies*

$$\lim_{t \to \infty} \sup_{0 \leq h \leq T} \|\phi_h X_t - X_{t+h}\| = 0$$

*for any $T > 0$.*
2. *The limit set $L(X_t)$ of $X_t$ is almost surely an internally chain recurrent set for the mean limit ODE.*

The first assertion of the theorem roughly states that $X_t$ tracks the flow of the mean limit ODE with increasing accuracy far into the future. The second assertion of the theorem states that the only candidates for limit sets of the process $x_n$ corresponding to the distribution of balls are connected compact internally chain recurrent sets for the mean limit flow.

To give a sense of the utility of this result, we derive some corollaries for the replicator processes and the fertility-selection process in the next two sections.

3.1. *Implications for replicator processes.* Let $\{R_n\}_{n\geq 0}$ and $\{\tilde{R}_n\}_{n\geq 0}$ be sequences of independent identically distributed random $k \times k$ matrices whose entries take values in the set $\{-1, 0, 1, \ldots, m\}$. Let $\{r_n\}_{n\geq 0}$ be a sequence of independent identically distributed random $k \times 1$ matrices whose entries take values in the set $\{-1, 0, 1, \ldots, m\}$. Let $z_n \in \mathbb{Z}_+^k$ be the replicator process associated with these random matrices. Define the mean payoff matrix by $A = E[R_0 + \tilde{R}_0]$. The limiting mean ODE associated with this process is given by a replicator equation [9]

$$\text{(6)} \qquad \frac{dx}{dt} = \text{diag}(x)Ax - (x.Ax)x, \qquad i = 1, \ldots, k,$$

where $x.A$ denotes multiplying the left-hand side of $A$ by the transpose of $x$ and $\text{diag}(x)$ is a diagonal $k \times k$ matrix with diagonal entries $x^i$. The dynamics of (6) are well studied and have two remarkable properties whose proofs can be found in [9].

THEOREM 2 (Exclusion principle). *If the replicator equation (6) has no equilibrium in* $\text{int } S_k$, *then every orbit of (6) converges to* $\partial S_k$.

THEOREM 3 (Time averaging principle). *If the replicator equation (6) has a unique equilibrium $p$ in* $\text{int } S_k$ *and if $x(t)$ is a solution of (6) such that $L(x(t)) \subset \text{int } S_k$, then*

$$\lim_{T \to \infty} \frac{1}{T} \int_0^T x(t)\, dt = p.$$

It turns out that Theorem 1 provides us the tool in which to transfer these theorems to replicator processes.

THEOREM 4. *Let $z_n$ be a replicator process on $\mathbb{Z}_+^k$ with mean payoff matrix $A$. If mean limit replicator equation has no equilibria in* $\text{int } S_k$, *then $L(x_n) \cap \partial S_k \neq \varnothing$ almost surely on the event $\{\liminf_{n\to\infty} \frac{|z_n|}{n} > 0\}$.*

PROOF. The proof of Theorem 2 implies there exists a vector $c \in \mathbb{R}^k$ such that the function $V(x) = \sum c^i \log x^i$ is strictly increasing along the forward orbits of the mean limit replicator equation that lie in $\text{int } S_k$. Consequently, every compact connected internally chain recurrent set intersects $\partial S_k$. Applying Theorem 1 completes the proof. □

THEOREM 5. *Let $z_n$ be a replicator process on $\mathbb{Z}_+^k$ with mean payoff matrix $A$. Let $X_t$ be continuous-time process associated with $z_n$ that is defined by (5). Assume that (6) has a unique rest point $p$ in* $\text{int } S_k$. *Then*

$$\lim_{T \to \infty} \frac{1}{T} \int_0^T X_t\, dt = p$$

almost surely on the event
$$\left\{\liminf_{n\to\infty}\frac{|z_n|}{n}>0\right\}\cap\{L(X_t)\subset\operatorname{int}S_k\}.$$

PROOF. Consider a trajectory $X_t$ from the event $\{\liminf_{n\to\infty}\frac{|z_n|}{n}>0\}\cap\{L(X_t)\subset\operatorname{int}S_k\}$. Theorem 1 implies that $X_t$ is almost surely an asymptotic pseudotrajectory for the flow $\phi_t$ of the mean limit replicator equation. Theorem 1 implies that $L(X_t)$ is a compact internally chain recurrent set for the flow $\phi_t$ of (6). Theorem 3 implies that
$$\lim_{t\to\infty}\frac{1}{t}\int_0^t\phi_s x\,ds=p$$
for all $x\in L(X_t)$, and this convergence is uniform. Therefore given $\varepsilon>0$, we can choose $T>0$ and a compact neighborhood $U$ of $L(X_t)$ such that
$$\left\|\frac{1}{T}\int_0^T\phi_s x\,ds-p\right\|<\frac{\varepsilon}{3},$$
whenever $x\in U$. Since $X_t$ is an asymptotic pseudotrajectory, there exists an $l\geq 1$ such that
$$\sup_{0\leq h\leq T}\|X_{t+h}-\phi_h X_t\|<\frac{\varepsilon}{3}$$
for all $t\geq lT$. For any $i\in\mathbb{Z}_+$, define
$$\psi(i)=\left\|\int_0^T(X_{iT+s}-\phi_s X_{iT})\,ds\right\|+\left\|\int_0^T(\phi_s X_{iT}-p)\,ds\right\|.$$
Since $L(X_t)\subset U$, there is an $N\geq l$ such that $X_t\in U$ for all $t\geq NT$. Given any $t\in\mathbb{R}$, let $[t]$ denote the integer part of $t$. For any $t>(N+1)T$, we get
$$\left\|\frac{1}{t}\int_0^t X_s\,ds-p\right\|$$
$$\leq\frac{1}{t}\left(\left\|\int_0^{NT}(X_s-p)\,ds\right\|+\psi(N)+\psi(N+1)+\cdots\right.$$
$$\left.+\psi\left(\left[\frac{t}{T}\right]-1\right)+\left\|\int_{[t/T]T}^t(X_s-p)\,ds\right\|\right)$$
$$\leq\frac{1}{t}\left(2NT+\varepsilon T\left(\left[\frac{t}{T}\right]-N\right)+2\left(t-\left[\frac{t}{T}\right]T\right)\right).$$
Taking the limit as $t\to\infty$, we get that
$$\limsup_{t\to\infty}\left\|\frac{1}{t}\int_0^t X_s\,ds-p\right\|\leq\varepsilon.$$
Taking the limit as $\varepsilon\to 0$ completes the proof of the theorem. $\square$

3.2. *Implications for additive fertility-selection processes.* Let $z_n \in \mathbb{Z}_+^{k \times k}$ be a fertility-selection process defined by the sequence of random variables $G_n(ij, rs)$ with $1 \leq i, j, r, s \leq k$. Define $g(ij, rs) = E[G_0(ij, rs)]$. Define $x_n = z_n/|z_n|$, whenever $z_n \neq 0$ and $x_n = 0$ otherwise.

The mean limit ODE for this selection-fertility process is given by the fertility-selection equations (see, e.g., [9])

$$\text{(7)} \qquad \frac{dx^{ij}}{dt} = \sum_{r,s=1}^{k} g(ir, js) x^{ir} x^{js} - x^{ij} \bar{g},$$

where

$$\bar{g} = \sum_{1 \leq i,j,r,s \leq k} g(ir, js) x^{ir} x^{js}.$$

Now consider the special case, when each allele contributes additively to number of progeny produced by a mating. In this case, if $\gamma_{ij}$ is genotype $A_i A_j$'s contribution to fertility, then the mating between genotypes $A_i A_j$ and $A_r A_s$ produces on average

$$g(ij, rs) = \gamma_{ij} + \gamma_{rs}$$

progeny. Under this additional assumption, equation (7) simplifies to

$$\text{(8)} \qquad \frac{dx^{ij}}{dt} = x^j \gamma_i + x^i \gamma_j - 2 x^{ij} \bar{\gamma},$$

where $x^i = \sum_{j=1}^{k} x^{ij}$ is the frequency of the allele $A_i$ in the population,

$$\gamma_i = \sum_{r=1}^{k} \gamma_{ir} x^{ir}$$

is the average fertility of allele $i$ in the population and

$$\bar{\gamma} = \sum_{i=1}^{k} \gamma_i = \sum_{i,j=1}^{k} \gamma_{ij} x^{ij}$$

is the average fertility of the population.

THEOREM 6. *If $z_n$ is an additive fertility-selection process and the mean-limit ODE* (8) *has only a finite number of equilibria, then on the event* $\{\liminf_{n \to \infty} \frac{|z_n|}{n} > 0\}$, $x_n$ *almost surely converges to an equilibrium of* (8).

The proof of this theorem follows from the work of Hofbauer and Sigmund [9] that we include here for the reader's convenience.

PROOF OF THEOREM 6. Define the Hardy–Weinberg manifold by

$$H = \{x : x^{ij} = x^i x^j \text{ for all } 1 \le i, j \le k\}.$$

Since for any solution $x(t)$ to (8)

$$\frac{d}{dt}(x^{ij}(t) - x^i(t)x^j(t)) = -(x^{ij}(t) - x^i(t)x^j(t))2\bar{\gamma},$$

$x^{ij}(t) - x^i(t)x^j(t)$ converges exponentially to zero. Hence, all compact connected internally chain recurrent sets lie in the Hardy–Weinberg manifold. On the Hardy–Weinberg manifold, the dynamics of (8) are determined by the Hardy–Weinberg relations $x^{ij} = x^i x^j$ and the differential equation

$$\frac{dx^i}{dt} = \gamma_i - x^i \bar{\gamma}. \tag{9}$$

This differential equation is the continuous-time selection equations with selection parameters $\gamma_{ij}$ and, consequently, the mean fertility $\bar{\gamma}$ is a strict Lyapunov function for (9) (see, e.g., [9]). Hence, all compact connected internally chain recurrent sets correspond to compact connected sets of equilibria. Since we have assumed that there are only a finite number of equilibria, the only compact connected internally chain recurrent sets are individual equilibria. Applying Theorem 1 completes the proof of this theorem. □

**4. Growth and convergence with positive probability.** Theorem 1 helps to determine the limiting behavior of the genotypic composition of a population on the event of growth. However, it does not indicate which limiting behaviors occur with positive probability and sheds no insight into conditions that ensure that the population grows with positive probability. The goal of this section is to show that when the mean limit ODE admits an attractor at which growth is expected, the population grows with positive probability and its genotypic composition converges to the attractor with positive probability. Prior to stating and proving this result, we prove the following proposition that estimates the rate of growth on the event of convergence to a set where growth is expected.

PROPOSITION 1. *Let $z_n$ be a generalized urn process satisfying assumptions* (A1) *and* (A2). *Let $K \subset S_k$ be a compact set. If*

$$\lambda = \inf_{x \in K} \sum_w p_w(x) \alpha(w) > 0,$$

*then*

$$\liminf_{n \to \infty} \frac{|z_n|}{n} \ge \lambda$$

*on the event $\{L(\{x_n\})_{n \ge 0} \subset K\} \cap \{\lim_{n \to \infty} |z_n| = \infty\}$.*

REMARK. If $K$ is an equilibrium, $\liminf_{n\to\infty} \frac{|z_n|}{n} = \lambda$ on the event $\{L(\{x_n\})_{n\geq 0} \subset K\} \cap \{\lim_{n\to\infty} |z_n| = \infty\}$.

PROOF OF PROPOSITION 1. Let

$$\mathcal{E} = \{L(\{x_n\})_{n\geq 0} \subset K\} \cap \left\{\lim_{n\to\infty} |z_n| = \infty\right\}.$$

We will show that $\liminf_{n\to\infty} |z_n|/n \geq \lambda - \varepsilon$ on the event $\mathcal{E}$ for every $\varepsilon > 0$. Let $\varepsilon > 0$ be given. The definition of $\lambda$, compactness of $K$, continuity of $p_w$ and assumption (A2) imply that there exist an integer $I$ and compact neighborhood $U$ of $K$ such that $E[|z_{n+1}| - |z_n||z_n = z] \geq \lambda - \varepsilon$, whenever $|z| \geq I$ and $z/|z| \in U$. For each natural number $j$, define the event $\mathcal{E}_j = \{|z_n| \geq I, z_n/|z_n| \in U \text{ for } n \geq j\}$. Notice that $\mathcal{E} \subset \bigcup_{j=1}^\infty \mathcal{E}_j$. Define a sequence of random variables by $N_0 = 0$ and

$$N_{n+1} = \begin{cases} |z_{n+1}| - |z_n|, & \text{if } |z_n| \geq I \text{ and } z_n/|z_n| \in U, \\ \lambda, & \text{otherwise}, \end{cases}$$

for $n \geq 0$. Let

$$M_n = \sum_{i=1}^n \frac{1}{i}(N_i - E[N_i|\mathcal{F}_{i-1}]).$$

$M_n$ is a martingale that satisfies

$$\sup_n E[M_n^2] \leq 4m^2 \sum_{i\geq 1} \frac{1}{i^2},$$

as $|N_i| \leq 2m$. Therefore, by Doob's convergence theorem $\{M_n\}_{n\geq 1}$ converges almost surely. By Kronecker's lemma,

(10) $$\lim_{n\to\infty} \frac{1}{n} \sum_{i=1}^n N_i - E[N_i|\mathcal{F}_{i-1}] = 0$$

almost surely. Since $\sum_{i=j+1}^n N_i = |z_n| - |z_j|$ for all $n \geq j$ on the event $\mathcal{E}_j$ and $E[N_i|\mathcal{F}_{i-1}] \geq \lambda - \varepsilon$ for all $i \geq 1$, (10) implies that $\liminf_{n\to\infty} \frac{|z(n)|}{n} \geq \lambda - \varepsilon$ almost surely on the event $\mathcal{E}_j$. It follows that $\liminf_{n\to\infty} \frac{|z(n)|}{n} \geq \lambda - \varepsilon$ almost surely on the event $\mathcal{E}$. □

Let $\phi_t(x)$ denote the flow of the mean limit ODE in (4). A compact invariant set $\mathcal{A} \subset S_k$ is called an *attractor* provided that there is an open neighborhood $U \subset S_k$ of $\mathcal{A}$ such that

$$\bigcap_{t>0} \overline{\bigcup_{s\geq t} \phi_s U} = \mathcal{A}.$$

The basin of attraction $B(\mathcal{A})$ of $\mathcal{A}$ is the set of points $x \in S_k$ satisfying $\inf_{y \in \mathcal{A}} \|\phi_t x - y\| \to 0$ as $t \to \infty$.

Define the set of *attainable points*, $\text{Att}_\infty(X)$, as the set of points $x \in S_k$ such that, for all $M \in \mathbb{N}$ and every open neighborhood $U$ of $x$

$$P[|z_n| \geq M \text{ and } x_n \in U \text{ for some } n] > 0.$$

THEOREM 7. *Let $z_n$ be a generalized urn process satisfying assumptions* (A1) *and* (A2). *Let $\mathcal{A}$ be an attractor for mean limit ODE with basin of attraction $B(\mathcal{A})$. Assume that*

$$\lambda = \inf_{x \in \mathcal{A}} \sum_w p_w(x) \alpha(w) > 0$$

*and define*

$$\mathcal{C} = \left\{ \liminf_{n \to \infty} \frac{|z_n|}{n} \geq \lambda \right\} \cap \{L(\{x_n\}_{n \geq 0}) \subseteq \mathcal{A}\}.$$

*If $U$ is an open set, whose closure is contained in $B(\mathcal{A})$, then there exists a constant $K > 0$ such that for all $M \in \mathbb{N}$,*

$$P[\mathcal{C}] \geq \left(1 - \frac{K}{M}\right) P[|z_n| \geq M \text{ and } x_n \in U \text{ for some } n].$$

*In particular, if*

$$B(\mathcal{A}) \cap \text{Att}_\infty(X) \neq \varnothing,$$

*then $P[\mathcal{C}] > 0$.*

REMARK. Theorem 7 simultaneously provides a condition that ensures that the population grows with positive probability and that the distribution of the population converges to an "attainable" attractor. Consequently, this result significantly improves ([12], Theorem 2.6) that proved convergence with positive probability to "attainable" attractor under the strong assumption of population growth with probability one.

PROOF OF THEOREM 7. Assumption $\inf_{x \in \mathcal{A}} \sum_w p_w(x) \alpha(w) > 0$ means that the population grows in a neighborhood of the attractor $\mathcal{A}$ when the population size is sufficiently large. It implies that there exist $a_1, a_2 > 0$ and a neighborhood $\mathcal{N}$ of $\mathcal{A}$ such that $E[|z_{n+1}| - |z_n||z_n] \geq a_1 \mathbb{1}_{\{x_n \in \mathcal{N}, |z_n| \geq a_2\}}$. The proof relies on the following principle: remaining in a neighborhood of the attractor increases the population size and this increase in population size increases the likelihood of remaining near the attractor. Let $U$ be an open set such that $\overline{U}$ is a compact subset of $B(\mathcal{A})$. Assume $M \in \mathbb{N}$ and $\rho \in \mathbb{N}$ are such that $P[|z_\rho| \geq M, x_\rho \in U] > 0$. Choose a neighborhood $V$ of $\mathcal{A}$ such that $\overline{V}$ is a compact subset of $\mathcal{N} \cap B(\mathcal{A})$. Since $\mathcal{A}$ is an attractor there exists a time $T_0 \geq 1$ such that the trajectories coming from $U \cup V$ rejoin the neighborhood $V$ after time $T_0$. More precisely, there exists $\delta > 0$ such that, if $X_t \in U \cup V$, $T \geq T_0$ and $\|\phi_T(X_t) - X_{t+T}\| < \delta$, then $X_{t+T} \in V$.

To avoid double subscripts, we let $z(r)$ denote $z_r$, $\tau(r)$ denote $\tau_r$, etc. Define $r_0 = \rho$ and

$$r_k = \inf\{r > r_{k-1}, \tau(r) - \tau(r_{k-1}) \geq T_0\}$$

for all $k \geq 1$. Since $\tau_{n+1} - \tau_n \geq \frac{1}{|z_0|+mn}$ for all $n$ where $m$ is the integer in assumption (A1), $\lim_{n\to\infty} \tau_n = \infty$ and $r_n < +\infty$ for all $n$. Define $A = e^{T_0 m}$, $B = 3A$ and the following events for all $k \geq 1$:

$E_1(k)$: $|z(r_k)| \geq \zeta^{k-1} B^{-1} M$,

$E_2(k)$: for all $r \in [r_k, r_{k+1}], x_r \in V$,

where $\zeta = 1 + a_1 T_0 / 2B$. Let $E(0)$ be the event $\{|z_\rho| \geq M, x_\rho \in U\}$. For $k \geq 1$, define $E(k) = E(k-1) \cap E_1(k) \cap E_2(k)$. We will show that there exists a constant $F > 0$ such that $P[E(k+1)|E(0)] \geq P[E(k)|E(0)] - F/M\zeta^k$ for all $k \geq 0$. The proof of this estimate relies on three lemmas. The first consists in observing that the population size $|z_r|$ remains bounded on time intervals of order $T_0$, namely between $B^{-1}|z(r_k)|$ and $B|z(r_k)|$ on $[r_k, r_{k+1}]$. The second one makes use of this claim to underestimate probability of being inside $V$ on $[r_{k+1}, r_{k+2}]$ if $x(r_k) \in U \cup V$. The third lemma estimates the probability that the population grows sufficiently.

LEMMA 2. *For enough large $|z(r_k)|$:*

1. $B|z(r_k)| \geq |z_r| \geq B^{-1}|z(r_k)|$ *for all $r \in [r_k, r_{k+1}]$,*
2. $T_0 \leq \tau(r_{k+1}) - \tau(r_k) \leq 2T_0$,
3. $T_0 B^{-1}|z(r_k)| \leq r_{k+1} - r_k \leq T_0 B |z(r_k)|$.

PROOF. Suppose $|z(r_k)| > Am$. Let $u = \inf\{n \in \mathbb{N} : |z(r_k + n)| \leq A^{-1}|z(r_k)|\}$. Then

$$\tau(r_k + u) - \tau(r_k) \geq \sum_{0 \leq j < |z(r_k)|(1-A^{-1})/m} \frac{1}{A^{-1}|z(r_k)| + mj}$$

$$\geq \int_0^{|z(r_k)|(1-A^{-1})/m} \frac{dx}{A^{-1}|z(r_k)| + mx}$$

$$= \frac{1}{m} \ln(A) = T_0.$$

This proves that $A^{-1}|z(r_k)| \leq |z(r)|$ for all $r \in [r_k, r_{k+1})$, which implies for sufficiently large $|z(r_k)|$ that $B^{-1}|z(r_k)| \leq |z(r)|$ for all $r \in [r_k, r_{k+1}]$.

One can show similarly that $2A|z(r_k)| \geq |z(r)|$ for all $r \in [r_k, r_{k+1}]$. The definition of $r_k$ immediately implies that $\tau(r_{k+1}) - \tau(r_k) \geq T_0$. Since $T_0 \geq 1$ and $\tau(n+1) - \tau(n) \leq 1$ for all $n$, the definition of $r_k$ also implies that $\tau(r_{k+1}) - \tau(r_k) \leq T_0 + 1 \leq 2T_0$. The proof of claim 1 and the fact that $T_0 \geq 1$ imply that $\frac{r_{k+1}-r_k}{|z(r_k)|B^{-1}} \geq \tau(r_{k+1}) - \tau(r_k) \geq T_0$ and $\frac{r-r_k}{2|z(r_k)|A} \leq \tau(r) - \tau(r_k) \leq T_0$ for all $r \in [r_k, r_{k+1})$. Claim 3 follows for sufficiently large $|z(r_k)|$. □

LEMMA 3. *There exists a $C > 0$ depending only on $p_w$, $T_0$, $\delta$, $a_1$ and $m$ such that*

$$P[E_2(k+1)|E_1(k), x(r_k) \in U \cup V] \geq 1 - \frac{C}{M\zeta^k}$$

*for all $k \geq 0$ and $M > 0$ sufficiently large.*

PROOF. Suppose $E_1(k)$ is satisfied and $x(r_k) \in U \cup V$. Lemma 2 implies that for large enough $M$, $|z(r)| \geq B^{-2}|z(r_k)|$ for $r \in [r_k, r_{k+2}]$, $4T_0 \geq \tau(r) - \tau(r_k) \geq T_0$ for $r \in [r_{k+1}, r_{k+2}]$ and $r_{k+2} - r_k \leq 2B^2 T_0 |z(r_k)|$. Define

$$g(x) = \sum_w p_w(x)(w - x\alpha(w)).$$

Let $L$ be the Lipschitz constant for $g$ and $\|g\|_0 = \sup \|g(x)\|$. Using Gronwall's inequality, we prove the following estimate:

$$(11) \quad \sup_{r \in [r_{k+1}, r_{k+2}]} \|\phi_{\tau(r) - \tau(r_k)} x(r_k) - x(r)\| \leq e^{4LT_0}(\Gamma_1(r_k, r_{k+2}) + \Gamma_2(r_k, r_{k+2})),$$

where

$$\Gamma_1(r_k, r_{k+2}) = \sup_{r_k \leq l \leq r_{k+2}-1} \left\| \sum_{i=r_k}^{l} \frac{U(i+1)}{|z(i)|} \right\|$$

and

$$\Gamma_2(r_k, r_{k+2}) = \frac{2\sup \|g(x)\|}{\inf_{r_k \leq r \leq r_{k+2}} |z(r)|} + \sup_{r_k \leq l \leq r_{k+2}-1} \left\| \sum_{i=r_k}^{l} \frac{b(i+1)}{|z(i)|} \right\|.$$

To prove (11), let $X(t) = X_t$ denote the continuous time version of $x_n$ defined in (5) and $c(t) = \sup\{n \in \mathbb{Z}_+ : t \geq \tau(n)\}$. Notice that for any $h \geq 0$ and $t \geq 0$,

$$X(t+h) - X(t)$$
$$= x(c(t+h)) - x(c(t)) = \sum_{i=c(t)}^{c(t+h)-1} x(i+1) - x(i)$$
$$= \sum_{i=c(t)}^{c(t+h)-1} \frac{g(x(i)) + U(i+1) + b(i+1)}{|z(i)|}$$

$$= \int_{\tau(c(t))}^{\tau(c(t+h))} g(X(s))\, ds + \sum_{i=c(t)}^{c(t+h)-1} \frac{U(i+1)+b(i+1)}{|z(i)|}$$

$$= \int_{\tau(c(t))-t}^{\tau(c(t+h))-t} g(X(t+s))\, ds + \sum_{i=c(t)}^{c(t+h)-1} \frac{U(i+1)+b(i+1)}{|z(i)|}.$$

Since $\phi_h X(t) = X(t) + \int_0^h g(\phi_s X(t))\, ds$, the previous equalities imply that

$$\|\phi_h(X(t)) - X(t+h)\|$$
$$\leq \int_0^h \|g(\phi_s(X(t))) - g(X(t+s))\|\, ds + \int_{\tau(c(t))-t}^0 \|g(X(t+s))\|\, ds$$
$$+ \int_{\tau(c(t+h))-t}^h \|g(X(t+s))\|\, ds + \left\|\sum_{i=c(t)}^{c(t+h)-1} \frac{U(i+1)+b(i+1)}{|z(i)|}\right\|$$
$$\leq L\int_0^h \|\phi_s(X(t)) - X(t+s)\|\, ds + \frac{\|g\|_0}{|z(c(t))|}$$
$$+ \frac{\|g\|_0}{|z(c(t+h))|} + \left\|\sum_{i=c(t)}^{c(t+h)-1} \frac{U(i+1)+b(i+1)}{|z(i)|}\right\|.$$

Choosing $t = \tau(r_k)$ and applying Gronwall's inequality to the previous inequality over the interval $0 \leq h \leq \tau(r_{k+2}) - \tau(r_k)$ gives the desired estimate.

Since $E_1(k)$ holds, $|z(r)| \geq |z(r_k)|B^{-2} \geq \zeta^{k-1}B^{-3}M$ for all $r \in [r_k, r_{k+2}]$. This observation, plus the fact that there exists $K > 0$ such that $\|b(n+1)\| \leq \frac{K}{|z(n)|}$, imply that $e^{4LT_0}\Gamma_2(r_k, r_{k+2}) < \frac{\delta}{2}$ for $M$ sufficiently large. On the other hand, Doob's inequality and Lemma 1 imply that

$$E\left[\sup_{r_k \leq l \leq r_{k+2}-1}\left\|\sum_{i=r_k}^l \frac{U(i+1)}{|z(i)|}\right\|^2 \bigg| z(r_k)\right] \leq 16m^2 E\left[\sum_{i=r_k}^{r_{k+2}-1} \frac{1}{|z(i)|^2}\bigg| z(r_k)\right]$$
$$\leq \frac{16m^2 B^4}{|z(r_k)|^2} E[r_{k+2} - r_k | z(r_k)]$$
$$\leq \frac{32m^2 T_0 B^6}{|z(r_k)|} \leq \frac{32m^2 B^7 T_0}{M\zeta^{k-1}}.$$

Therefore,

$$(12) \quad P\left[\sup_{r_k \leq l \leq r_{k+2}-1}\left\|\sum_{i=r_k}^l \frac{U(i+1)}{|z(i)|}\right\| \geq e^{-4LT_0}\frac{\delta}{2}\bigg| z(r_k)\right] \leq \frac{128m^2 T_0 B^7 e^{8LT_0}}{\delta^2 M\zeta^{k-1}}.$$

Define
$$\mathcal{E} = \left\{ \sup_{r \in [r_{k+1}, r_{k+2}]} d(\phi_{\tau(r) - \tau(r_k)} x(r_k), x(r)) \leq \delta \right\}.$$

Since $\tau(r_{k+1}) - \tau(r_k) \geq T_0$ and $x(r_k) \in U \cup V$, our choice of $T_0$ implies that $x(r) \in V$ for all $r \in [r_{k+1}, r_{k+2}]$ on the event $\mathcal{E}$. Inequalities (11) and (12) imply that $P[\mathcal{E}] \geq 1 - \frac{128 m^2 T_0 B^7 e^{8LT_0}}{\delta^2 M \zeta^{k-1}}$ for $M$ sufficiently large. □

LEMMA 4. *There exists $D > 0$, depending only on $p_w$, $T_0$, $\delta$, $a_1$ and $m$ such that*
$$P[E_1(k+1) \cup E_2(k)^c | E_1(k)] \geq 1 - \frac{D}{M\zeta^k}$$

*for all $k \geq 1$ and $M > 0$ sufficiently large.*

PROOF. Define
$$N(i) = |z(i)| - |z(i-1)|, \qquad D(i) = N(i) - E[N(i)|z(i-1)]$$

and
$$G(k+1) = \frac{1}{r_{k+1} - r_k} \sum_{i=r_k+1}^{r_{k+1}} D(i).$$

Observe that $|D(n)| \leq 2m$. Therefore,

$$E[G(k+1)^2 | z(r_k)] \leq \frac{B^2}{T_0^2 |z(r_k)|^2} E\left[\left(\sum_{i=r_k+1}^{r_k + T_0 B |z(r_k)|} D(i) \mathbb{1}_{\{i-1 < r_{k+1}\}}\right)^2 \bigg| z(r_k)\right]$$

$$\leq \frac{B^2}{T_0^2 |z(r_k)|^2} \sum_{i=r_k+1}^{r_k + T_0 B |z(r_k)|} E[D(i)^2 | z(r_k)] \leq \frac{4m^2 B^3}{T_0 |z(r_k)|},$$

where we have used the fact that $T_0 B^{-1} |z(r_k)| \leq r_{k+1} - r_k \leq T_0 B |z(r_k)|$. It follows that $P[G(k+1) \leq -\frac{a_1}{2}] \leq P[G(k+1)^2 \geq \frac{a_1^2}{4}] \leq \frac{D}{M\zeta^k}$ with $D = \frac{16 m^2 B^4 \zeta}{a_1^2 T_0}$. Since $\zeta = 1 + a_1 T_0 / 2B$, it follows that
$$E_1(k+1)^c \cap E_1(k) \cap E_2(k) \subset \{G(k+1) \leq -a_1/2\}. \qquad \square$$

These three lemmas imply that there exists $F > 0$ depending only on $p_w$, $T_0$, $\delta$, $a_1$ and $m$ such that
$$P[E(k+1)|E(0)] \geq P[E(k)|E(0)] - \frac{F}{M\zeta^k}$$

for all $k \geq 0$ and $M > 0$ sufficiently large. Indeed, due to the fact that for all $k \geq 1$ $E(k)$ equals the disjoint union of $E(k+1)$, $E_1(k+1)^c \cap E(k)$, and $E_2(k+1)^c \cap E_1(k+1) \cap E(k)$, we get

$$\begin{aligned}
P[E(k+1)|E(0)] \\
&= P[E(k)|E(0)] - P[E_2(k+1)^c \cap E_1(k+1) \cap E(k)|E(0)] \\
&\quad - P[E_1(k+1)^c \cap E(k)|E(0)] \\
&\geq P[E(k)|E(0)] \\
&\quad - P[E_2(k+1)^c \cap E_1(k) \cap \{x(r_k) \in U \cup V\}|E(0)] \\
&\quad - P[E_1(k+1)^c \cap E_2(k)|E(0)] \geq P[E(k)|E(0)] - \frac{F}{M\zeta^k},
\end{aligned}$$

where the first inequality follows from the inclusions $E(k) \cap E_1(k+1) \subset E_1(k) \cap \{x(r_k) \in U \cup V\}$ and $E(k) \subset E_2(k)$, and the second equality follows from Lemmas 3 and 4 with $F = C + D$. These inequalities remain true for $k = 0$, since $E_1(1)$ always holds.

It follows that

$$P\left[\lim_{k \to \infty} E(k)\right] \geq P[E(0)]\left(1 - \sum_{k=0}^{\infty} \frac{F}{M\zeta^k}\right) \geq P[E(0)]\left(1 - \frac{\zeta F}{M(\zeta - 1)}\right).$$

The definition of $E(k)$ implies that

$$P[\mathcal{C}] \geq P[E(0)]\left(1 - \frac{\zeta F}{M(\zeta - 1)}\right),$$

where $\mathcal{C} = \{\liminf_{n \to \infty} \frac{|z_n|}{n} > 0\} \cap \{x_n \in U \cup V \text{ i.o.}\}$. On the event $\mathcal{C}$, Theorem 1 implies that $L(\{x_n\})$ is a compact internally chain recurrent set for the mean limit ODE. Since $L(\{x_n\}) \cap B(\mathcal{A}) \neq \varnothing$ on the event $\mathcal{C}$, a basic result about internally chain recurrent sets (see, e.g., [4], Corollary 5.4) implies that $L(\{x_n\}) \subset \mathcal{A}$ on the event $\mathcal{C}$. Setting $K = \frac{\zeta F}{M(\zeta - 1)}$ and applying Proposition 1 completes the proof of the first assertion of the theorem. To prove the second assertion, assume that $p \in \text{Att}_\infty(X) \cap B(\mathcal{A})$, choose $U$ an open neighborhood of $p$ such that $\bar{U} \subset B(\mathcal{A})$, and apply the first assertion of the theorem. $\square$

**5. Nonconvergence.** In this section we show that there are two types of invariant sets of the mean limit ODE toward which the generalized urn process does not converge. The first type corresponds to a compact set where growth of the process is not expected, and the second type corresponds to a "nondegenerate" equilibrium or periodic orbit.

PROPOSITION 2. *Let $z_n$ be a Markov process on $\mathbb{Z}_+^k$ satisfying assumptions* (A1) *and* (A2). *If $K \subset S_k$ is a compact set satisfying*

$$\sup_{x \in K} \sum_w p_w(x) \alpha(w) < 0, \tag{13}$$

*then $P[\{L(x_n) \subset K\} \cap \{\lim_{n \to \infty} |z_n| = +\infty\}] = 0$.*

PROOF. Equation (13) implies that we can choose a neighborhood $U$ of $K$, $N \in \mathbb{N}$ and $\varepsilon > 0$ such that

$$\sup_{|z| \geq N, z/|z| \in U} E[|z_{n+1}| - |z_n| | z_n = z] \leq -\varepsilon. \tag{14}$$

Given any $l \in \mathbb{N}$ such that $x_l \in U$ and $|z_l| \geq N$, define the stopping time

$$T = \inf\{n \geq l : x_n \notin U \text{ or } |z_l| < N\}.$$

For any $n \geq l$, we get

$$0 \leq E[|z_{T \wedge n}|] = E\left[\sum_{i=l+1}^{T \wedge n} |z_i| - |z_{i-1}|\right] + E[|z_l|]$$

$$= \sum_{i=l+1}^n E[|z_{i \wedge T}| - |z_{(i-1) \wedge T}|] + E[|z_l|]$$

$$\leq -\varepsilon \sum_{i=l+1}^n P[T \geq i] + E[|z_l|].$$

Taking the limit as $n \to \infty$, we get that $\sum_{i=l+1}^\infty P[T \geq i] \leq E[|z_l|]/\varepsilon$. The Borel–Cantelli lemma implies that $P[T = \infty] = 0$. It follows that $P[\{L(x_n) \subset K\} \cap \{\lim_{n \to \infty} |z_n| = +\infty\}] = 0$. □

If $A$ is a subset of $\mathbb{R}^k$, then we let $\text{Span}(A) \subset \mathbb{R}^k$ denote the vector space spanned by the points in $A$. Given a compact subset $\mathcal{U} \subset \text{int}(S_k)$, we say that the process $\{z_n\}$ is *nondegenerate at $\mathcal{U}$* if for all $x \in \mathcal{U}$,

$$\text{Span}\{w \in \mathbb{Z}^k : p_w(x) > 0\} = \mathbb{R}^k.$$

Recall, a periodic orbit or an equilibrium of an ODE is *linearly unstable* provided that one of its characteristic exponents is greater than zero.

THEOREM 8. *Let $\{z_n\}$ be a generalized urn process satisfying assumptions* (A1) *and* (A2). *Let $\mathcal{U} \subset \text{int}(S_k)$ be a linearly unstable periodic orbit or equilibrium for the mean limit ODE. Assume the following:*

(a) *There exists $\beta > 1/2$ such that the functions $p_w$ are $C^{1+\beta}$ in a neighborhood of $\mathcal{U}$.*

(b) $\{z_n\}$ is nondegenerate at $\mathcal{U}$.

Then $P[(L(x_n) \subset \mathcal{U}) \cap \{\liminf_{n\to\infty} \frac{|z_n|}{n} > 0\}] = 0$.

PROOF. Let $N(\mathcal{U})$ be a neighborhood of $\mathcal{U}$. The event

$$\left\{L(x_n) \subset \mathcal{U} \text{ and } \liminf_{n\to\infty} \frac{|z_n|}{n} > 0\right\}$$

is contained in the event

$$\bigcup_{N\in\mathbb{Z}_+, \lambda\in\mathbb{Q}_+^*} E_{N,\lambda},$$

where $\mathbb{Q}_+^*$ denotes the positive rationales and

$$E_{N,\lambda} = \{L(x_n) \subset \mathcal{U}\} \cap \{\forall n \geq N \, |z_n| \geq n\lambda \text{ and } x_n \in N(\mathcal{U})\}.$$

In order to prove Theorem 8 it then suffices to prove that for $N$ large enough and $\lambda \in \mathbb{Q}_+^*$,

$$P[E_{N,\lambda}] = 0.$$

Let $\mathcal{F}_n$ denote the sigma field generated by $z_0, \ldots, z_n$ and $V_{n+1} = |z_n|(x_{n+1} - x_n)$. Let $F$ denote the vector field on $S_k$ defined by $F(x) = \sum_w p_w(x)(w - x\alpha(w))$. Let $\{\varepsilon_n\}$ denote a sequence of bounded, zero-mean i.i.d. random variables taking values in

$$TS_k = \left\{u \in \mathbb{R}^k : \sum u_i = 0\right\},$$

whose covariance matrix is nondegenerate (i.e., has rank $k-1$).

Define the sequence $\{\tilde{x}_n\}_{n\geq N}$ as follows:

$$\tilde{x}_N = x_N,$$

(15) $$\tilde{x}_{n+1} - \tilde{x}_n = \begin{cases} \frac{1}{|z_n|}(F(\tilde{x}_n) + V_{n+1} - F(x_n)), \\ \qquad\qquad\qquad\qquad \text{if } x_n \in N(\mathcal{U}) \text{ and } |z_n| \geq n\lambda, \\ \frac{1}{\lambda n}(F(\tilde{x}_n) + \varepsilon_{n+1}), \qquad \text{otherwise.} \end{cases}$$

The processes $\{x_n\}$ and $\{\tilde{x}_n\}$ coincide on the event $E_{N,\lambda}$. On the other hand,

$$P\left[\lim_{n\to\infty} \text{dist}(\tilde{x}_n, \mathcal{U}) = 0\right] = 0$$

in view of the following theorem whose proof is an easy adaptation of [13], Theorem 2.

THEOREM 9. *Let $\{\mathcal{F}_n\}_{n\in\mathbb{Z}^+}$ denote a nondecreasing sequence of sub-$\sigma$-algebras of $\mathcal{F}$ and $(\tilde{x}_n)$ a sequence of adapted random variables given by*

(16) $$\tilde{x}_{n+1} - \tilde{x}_n = \beta_n(F(\tilde{x}_n) + \tilde{U}_{n+1} + \tilde{b}_{n+1}),$$

*where $F$ is a $C^{1+\beta}$ vector field, with $1/2 < \beta \leq 1$, $\{\tilde{U}_n\}$, $\{\tilde{b}_n\}$ and $\{\beta_n\}$ are adapted random variables.*

*We define $\tilde{V}_{n+1} = F(\tilde{x}_n) + \tilde{U}_{n+1} + \tilde{b}_{n+1}$.*

*Assume the following:*

(i) *$\exists K > 0$, $\forall n \in \mathbb{Z}_+$ $\|\tilde{U}_n\| \leq K$ and $E(\tilde{U}_{n+1}|\mathcal{F}_n) = 0$.*
(ii) *There exist $a$, $b > 0$ and a deterministic sequence $\{\gamma_n\}$ of nonnegative numbers having infinitely positive terms, such that $\forall n \in \mathbb{Z}_+, a\gamma_n \leq \beta_n \leq b\gamma_n$.*
(iii) *$\sum \tilde{b}_i^2 < +\infty$.*
(iv) *$\mathcal{U} \subset \mathrm{Int}(S_k)$ is a linearly unstable periodic orbit or equilibrium for $F$.*
(v) *There exist a neighborhood $N(\mathcal{U})$ of $\mathcal{U}$ and $c > 0$ such that, for all unit vectors $v \in \mathbb{R}^m$, $E(|\langle \tilde{V}_{n+1}, v \rangle\| \mathcal{F}_n) \geq c\mathbb{1}_{\{\tilde{x}_n \in N(\mathcal{U})\}}$.*

*Then $P[L(\tilde{x}_n) \subset \mathcal{U}] = 0$.*

Let
$$A_n = \{x_n \in N(\mathcal{U}) \text{ and } |z_n| \geq n\lambda\}.$$
Using the notation of Lemma 1 set
$$\tilde{U}_{n+1} = U_{n+1}, \qquad \tilde{b}_{n+1} = b_{n+1}, \qquad \beta_n = 1/|z_n| \qquad \text{on } A_n$$
and
$$\tilde{U}_{n+1} = \varepsilon_{n+1}, \qquad \tilde{b}_{n+1} = 0, \qquad \beta_n = 1/\lambda n \qquad \text{on } A_n^c.$$
Then, by Lemma 1, the process $\{\tilde{x}_n\}$ defined by (15) verifies recursion (16) and assertions (i)–(iv) of Theorem 9 are satisfied.

It remains to verify assertion (v). Let $B(1) = \{v \in S_k : \|v\| = 1\}$. Let $G_1 : S_k \times B(1) \to \mathbb{R}_+$ and $G_2 : \mathcal{U} \times S_k \times B(1) \to \mathbb{R}_+$ be the functions defined by
$$G_1(\tilde{x}, v) = E(|\langle F(\tilde{x}) + \varepsilon_n, v \rangle|)$$
and
$$G_2(x, \tilde{x}, v) = \sum_e |\langle F(\tilde{x}) - F(x) + Q_x(e), v \rangle| p_e(x),$$
where $Q_x$ denote the projection operator $Q_x : \mathrm{Span}\{x\} \oplus TS_k \to TS_k$.

Then it is not hard to verify that
$$E(|\langle \tilde{V}_{n+1}, v \rangle\| \mathcal{F}_n) = G_1(\tilde{x}_n, v)\mathbf{1}_{A_n^c} + G_2(x, \tilde{x}_n, v)\mathbf{1}_{A_n} + O(1/n).$$
By continuity of $G_1$, $G_2$ compactness of the sets $\mathcal{U}, S_k, B(1)$ and assumption (A2), there exists $b > 0$ and a neighborhood $N(\mathcal{U})$ such that $G_1(\tilde{x}, v) > b$, $G_2(x, \tilde{x}, v) > b$ for all $x \in N(\mathcal{U})$, $\tilde{x} \in S_k$ and $v \in B(1)$

Assumption (v) is thus verified for $n \geq N$ and $N$ large enough. $\square$

**6. Nondegenerate processes with gradient-like mean limit ODEs.** Using the results from the previous two sections, we can prove the following result.

THEOREM 10. *Let $z_n$ be a generalized urn process satisfying assumptions (A1) and (A2), $p_w$ are $C^{1+\beta}$ for some $\beta > 1/2$, $z_n$ is nondegenerate on $S_k$, $\text{Att}_\infty(X) = S_k$ and the chain recurrent set for the mean limit ODE consists of hyperbolic equilibria $q$ that satisfy $\sum_w p_w(q)\alpha(w) \neq 0$. Then $P[\liminf_{n\to\infty} \frac{|z_n|}{n} > 0]$ if and only if there exists a linearly stable equilibrium $q$ such that $\sum_w p_w(q)\alpha(w) > 0$. Furthermore, $P[\{L(x_n) = q\}] > 0$ for any linearly stable equilibrium $q$ satisfying $\sum_w p_w(q)\alpha(w) > 0$, and $P[\{L(x_n) = q\}] = 0$ for any equilibrium $q$ which is linearly unstable.*

PROOF. Define $\mathcal{G} = \{\liminf_{n\to\infty} |z_n|/n > 0\}$. Since $z_n$ is nondegenerate, Theorems 1 and 8 imply that $L(x_n)$ is contained almost surely in the set of linearly stable equilibria on the event $\mathcal{G}$. If $\sum_w p_w(q)\alpha(w) < 0$ for all linearly stable equilibria $q$, then Proposition 2 implies $P[\mathcal{G}] = 0$. Alternatively, if $q$ is a linearly stable equilibrium and $\sum_w p_w(q)\alpha(w) > 0$, then Theorem 7 implies $P[\mathcal{G} \cap \{L(x_n) = q\}] > 0$. Finally, if $q$ is an equilibrium which is linearly unstable, then Proposition 2 implies $P[L(x_n) = q] = 0$ if $\sum_w p_w(q)\alpha(w) < 0$ and Proposition 1 and Theorem 8 imply $P[L(x_n) = q] = 0$ if $\sum_w p_w(q)\alpha(w) > 0$. □

As an application of this result, we consider additive fertility-selection processes $z_n$ with mutation, where $g(ij, rs) = E[G_n(ij, rs)]$ is the expected number of progeny produced by a mating between genotypes $A_i A_j$ and $A_r A_s$, and $\mu(ij, rs)$ is the probability genotype $A_i A_j$ mutates to genotype $A_r A_s$.

COROLLARY 1. *Let $z_n$ be an additive fertility-selection process with mutation. Suppose:*

- *$\forall i, j, r, s \in \{1, \ldots, k\}$, $\mu(rs, ij)$ is strictly positive, and sufficiently small when $\{r, s\} \neq \{i, j\}$.*
- *$\forall i, j, r, \ s \in \{1, \ldots, k\}$ and such that the fertility-selection equation (8) without mutation has hyperbolic equilibria $q$ satisfying $\sum_w p_w(q)\alpha(w) \neq 0$, $P[G_1(ij, rs) \geq 3] > 0$.*

*Then $P[\liminf_{n\to\infty} \frac{|z_n|}{n} > 0]$, if and only if there exists a linearly stable equilibrium $q$ such that $\sum_w p_w(q)\alpha(w) > 0$. Furthermore, $P[\{L(\{x_n\}) = q\}] > 0$ for any linearly stable equilibrium $q$ satisfying $\sum_w p_w(q)\alpha(w) > 0$, and $P[\{L\{x_n\} = q\}] = 0$ for any equilibrium $q$ which is linearly unstable.*

Proof. Since the mean limit ODE corresponding to the fertility-selection process without mutation is gradient-like and has only hyperbolic equilibria, the chain-recurrent set for this mean limit ODE equals the set of equilibria. Consequently, the mean limit ODE corresponding to the fertility-selection process with sufficiently small mutation rates also has a chain-recurrent set consisting only hyperbolic equilibria. Due to the fact that all mutation rates are positive, this process is nondegenerate on the entire simplex. Since $P[G_1(ij,rs) \geq 3] > 0$ for all $1 \leq i,j,r,s \leq k$, $\text{Att}_\infty(X)$ is the entire simplex. Applying Theorem 10 completes the proof. $\square$

M. Benaïm
Faculté des sciences
Institut de Mathématiques
Rue Emile-Argand 11
Case Postale 2
CH-2007 Neuchâtel Suisse
France
e-mail: michel.benaim@unine.ch

S. J. Schreiber
Department of Mathematics
College of William and Mary
Williamsburg, Virginia 23187-8795
USA
e-mail: sjs@math.wm.edu

P. Tarrès
Laboratoire de Statistique
  et Probabilités
UMR C5583–Université Paul Sabatier
Bât 1R1-118, route de Narbonne
31062 Toulouse Cedex 4
France
e-mail: tarres@math.ups-tlse.fr